\documentclass[12pt]{article}

\usepackage[geometry]{ifsym}
\usepackage{verbatim}
\usepackage{fancyhdr}
\usepackage{graphicx}
\usepackage{mathrsfs}
\usepackage{amssymb}
\usepackage{cite}
\usepackage{epsfig}
\usepackage{pst-poly}     
\usepackage{pst-all}      

\usepackage{amsfonts,amssymb, amsmath, latexsym}

\newtheorem{pps}{Proposition}[section]
\newtheorem{cor}{Corollary}[section]
\newtheorem{lem}{Lemma}[section]
\newtheorem{thm}{Theorem}[section]

\newenvironment{pf}[1][Proof]{\noindent\textbf{#1.} }{\hfill\rule{1mm}{2mm}}

\makeatletter \@addtoreset{equation}{section} \makeatother

\begin{document}

\title{On the diameter of the Kronecker product graph\thanks {The work was supported by NNSF
of China (No. 11071233).}}
\author
{Fu-Tao Hu,\quad Jun-Ming Xu\footnote{Corresponding author:
xujm@ustc.edu.cn}\ \\
{\small School of Mathematical Sciences}  \\
{\small University of Science and Technology of China}\\
{\small Wentsun Wu Key Laboratory of CAS}\\
{\small Hefei, Anhui, 230026, China} }
\date{}
\maketitle

\begin{quotation}
\textbf{Abstract}: Let $G_1$ and $G_2$ be two undirected nontrivial
graphs. The Kronecker product of $G_1$ and $G_2$ denoted by
$G_1\otimes G_2$ with vertex set $V(G_1)\times V(G_2)$, two vertices
$x_1x_2$ and $y_1y_2$ are adjacent if and only if $(x_1,y_1)\in
E(G_1)$ and $(x_2,y_2)\in E(G_2)$. This paper presents a formula for
computing the diameter of $G_1\otimes G_2$ by means of the diameters
and primitive exponents of factor graphs.

\vskip6pt\noindent{\bf Keywords}: diameter, Kronecker product,
primitive exponent

\noindent{\bf AMS Subject Classification: }\ 05C12, 05C50

\end{quotation}

\section{Introduction}

For notation and graph-theoretical terminology not defined here we
follow \cite{x03}. Specifically, let $G=(V,E)$ be a nontrivial graph
with no parallel edges, but loops allowed, where $V=V(G)$ is the
vertex-set and $E=E(G)$ is the edge-set.

For two graphs $G$ and $H$, Kronecker product $G\otimes H$ is a
graph with vertex set $V(G)\times V(H)$ and two vertices $x_1x_2$
and $y_1y_2$ are adjacent when $(x_1,y_1)\in E(G)$ and $(x_2,y_2)\in
E(G)$.

As an operation of graphs, Kronecker product $G\otimes H$ was
introduced first by Weichesel~\cite{w62} in 1962. It has been shown
that the Kronecker product is a good method to construct lager
networks that can generate many good properties of the factor graphs
(see~\cite{lckfc10}), and has received much research attention
recently. Some properties and graphic parameters have been
investigated~\cite{bm98, bkr07,egk08,kz96,mv08}. The connectivity
and diameter are two important parameters to measure reliability and
efficiency of a network. Very recently, the connectivity of
Kronecker product graph has been deeply studied (see,
\cite{bs08,gqg10,gv09,mv08,o11,ww11,wx11,wy12}). However, the
diameter of Kronecker product graph has been not investigated yet.

In this paper, we determine the diameter of Kronecker product graph
by means of primitive exponents and diameters of factor graphs. In
particular, we obtain that
 $$
 d(G_1\otimes G_2)=\left\{ \begin{array}{ll}
\gamma_1\ & {\rm if }\ \gamma_1=\gamma_2; \\
\max\{\gamma_2+1,d_1\}\ & {\rm if }\ \gamma_1>\gamma_2;\\
\max\{\gamma_1+1,d_2\}\ & {\rm if }\ \gamma_1<\gamma_2,
\end{array}
 \right.
 $$
where $\gamma_i$ and $d_i$ are the primitive exponent and diameter
of $G_i$ for $i=1,2$, respectively.


\section{Some Lemmas}

Let $G$ be a graph. Denote $\gamma(G;x,y)$ to be the minimum integer
such that there exists an $(x,y)$-walk of length $k$ for any $k\ge
\gamma(G;x,y)$ and $\gamma(G)$ be the minimum integer $\gamma$ for
which, for any two vertices $x$ and $y$ in $G$, there exists an
$(x,y)$-walk of length $k$ for any integer $k\ge \gamma$. Let
$\gamma(G)=\max\{\gamma(G;x,y): x,y\in V(G)\}$.

If $\gamma(G)$ is well-defined, then $G$ is said to be {\it
primitive}, and $\gamma(G)$ is called the {\it primitive exponent},
{\it exponent} for short, of $G$. If $\gamma(G)$ does not exist,
then denote $\gamma(G)=\infty$.

Let $K_n^+$ be a graph obtained from a complete graph $K_n$ by
appending a loop on each vertex. It is clear that for a graph $G$
without parallel edges of order $n$, $\gamma(G)=1$ if and only if
$G\cong K_n^+$.

Let $A$ be the adjacency matrix of $G$. Equivalently, the exponent
of $G$ is the minimum integer $\gamma$ for which $A^\gamma>0$ and
$A^k\ngtr 0$ for any positive integer $k<\gamma$. Let $A_i$ be  the
adjacent matrix of $G_i$ for $i=1,2$. Since for any positive integer
$k$, $(A_1\otimes A_2)^k=A_1^k\otimes A_2^k$, by definition, we have
the following result immediately.

\begin{pps}\label{lem2.9a}   
Let $G_i$ be a primitive graph with exponent $\gamma_i$ for $i=1,2$,
and $G=G_1\otimes G_2$. Then $\gamma(G)=\max\{\gamma_1,\gamma_2\}$.
\end{pps}


The following lemmas will be used in proofs of our main results.

\begin{lem}\label{lem2.1} {\rm (Liu et al.~\cite{lmwz90})}
A graph $G$ is primitive if and only if $G$ is connected
and contains odd cycles.
\end{lem}

\begin{lem}\label{lem2.2} {\rm (Liu et al.~\cite{lmwz90})}
Let $G$ be a primitive graph, and let $x$ and $y$ be any pair of
vertices in $V(G)$. If there are two $(x,y)$-walks $P_1$ and $P_2$
with lengths $k_1$ and $k_2$, respectively, where $k_1$ and $k_2$
have different parity, then $\gamma(G;x,y)\le \max\{k_1, k_2\}-1$.
\end{lem}

\begin{lem}\label{lem2.3} {\rm (Delorme and Sol\'{e}~\cite{ds91})}
If $G$ is a primitive graph with diameter $d$, then $\gamma(G)\le
2d$.
\end{lem}

\begin{lem}\label{lem2.4} {\rm (Weichesel~\cite{w62})}
Let $G_1$ and $G_2$ be two connected graphs
and $G=G_1\otimes G_2$. Then $G$ is connected
if and only if either $G_1$ or $G_2$ contains an odd cycle.
\end{lem}

\begin{lem}\label{lem2.5}
Let $G=G_1\otimes G_2$, $x_i$ and $y_i$ be any two vertices in $G_i$,
$P_i$ be an $(x_i,y_i)$-walk of length $\ell_i$ in $G_i$
for $i=1,2$. If $\ell_1$ and $\ell_2$ have same parity, then
there is an $(x_1x_2,y_1y_2)$-walk of length
$\max\{\ell_1,\ell_2\}$ in $G$.
\end{lem}

\begin{pf}
Without loss of generality, suppose $\ell_1\ge \ell_2$.
Let $k=\ell_1-\ell_2$. Then $k$ is even.
Let $P_1=(x_1,z_1\ldots,z_{\ell_1-1},y_1)$ and
$P'_2=(x_2,u_1\ldots,u_{\ell_1-1},y_2)$ be an $(x_2,y_2)$-walk of length
$\ell_1$ in $G_2$ obtained from $P_2$ by repeating $k$ times of some edge in $P_2$.
Then
 $
 (x_1x_2, z_1u_1,\ldots,z_{\ell_1-1}u_{\ell_1-1},y_1y_2)
 $
is an $(x_1x_2,y_1y_2)$-walk of length $\ell_1$ in $G$.
\end{pf}

\begin{lem}\label{lem2.6}
Let $G$ be a primitive graph with exponent $\gamma$ and
order $n\ge 2$. We have\\
(i) if $\gamma$ is odd, then there exist two vertices $x$ and $y$,
and two different vertices $u$ and $v$, such that
the shortest odd $(x,y)$-walk and the shortest even $(u,v)$-walk are of
length $\gamma$ and $\gamma+1$, respectively;\\
(ii) if $\gamma$ is even, then there exist two different vertices $p$ and $q$,
and two vertices $w$ and $s$, such that
the shortest even $(p,q)$-walk and the shortest odd $(w,s)$-walk are of
length $\gamma$ and $\gamma+1$, respectively.
\end{lem}

\begin{pf}
(i) Assume that $\gamma$ is odd. If $\gamma=1$, then $G\cong K_n^+$.
Let $u$ and $v$ be two different vertices in $G$. Then the shortest
odd $(u,v)$-walk and the shortest even $(u,v)$-walk are of length 1
and $2$, respectively. Suppose now $\gamma\ge 3$.

Let $A$ be the adjacency matrix of $G$. By definition of $\gamma$,
$A^{\gamma-1}\ngtr 0$ and $A^{\gamma-2}\ngtr 0$. These imply that
there exist four vertices $x$, $y$, $u$ and $v$ such that there are
no odd $(x,y)$-walk and even $(u,v)$-walk with length $\gamma-2$ and
$\gamma-1$, respectively. Hence there are no odd $(x,y)$-walk and
even $(u,v)$-walk with length no more than $\gamma-2$ and
$\gamma-1$, respectively. Therefore, the shortest odd $(x,y)$-walk
and the shortest even $(u,v)$-walk are of length $\gamma$ and
$\gamma+1$, respectively.

We now show $u\ne v$.
If $u=v$, then $(u,w,u)$ is an even $(u,v)$-walk of length $2$
for any vertex $w$ adjacent to $u$ in $G$,
a contradiction with $\gamma+1\ge 3$.

(ii) Assume $\gamma$ is even. If $\gamma=2$, then $d=d(G)=1$ or $2$
since $d\le \gamma$. If $d=1$, then $G$ is isomorphic to a complete
graph $K_n$ with $m$ vertices having loops and $m<n$. Let $p$ be a
vertex with no loop and $q\ne p$ be another vertex in $G$. Then the
shortest even $(p,q)$-walk and odd $(p,p)$-walk are of length $2$
and $3$, respectively. If $d=2$, then there exist two different
vertices $p$ and $q$ such that $d_G(p,q)=2$, and hence the shortest
even $(p,q)$-walk and odd $(p,q)$-walk are of length $2$ and $3$,
respectively.

The case when $\gamma>2$ can be proved
by applying the similar discussion as in (i).
\end{pf}

\begin{lem}\label{lem2.7}
Let $G_i$ be a primitive graph with exponent $\gamma_i$ for $i=1,2$,
$G=G_1\otimes G_2$, and $x=x_1x_2$ and $y=y_1y_2$ be two different
vertices in $G$. If the shortest odd
(resp. even) $(x_1,y_1)$-walk in $G_1$ and the shortest even (resp. odd)
$(x_2,y_2)$-walk in $G_2$ are of length $m$ and $n$, respectively,
then $d_G(x,y)\ge \min\{m,n\}$.
\end{lem}
\begin{pf}
Without loss of generality, assume that $m$ is odd and $n$ is even.
Let
$$P=(x_1x_2,\ldots,u_1u_2,\ldots,y_1y_2)$$
be a minimum $(x,y)$-path with length $s$ in $G$.
Then
$$(x_1,\ldots,u_1,\ldots,y_1)\ {\rm and} \ (x_2,\ldots,u_2,\ldots,y_2)$$
be an $(x_1,y_1)$-walk in $G_1$ and an $(x_2,y_2)$-walk in $G_2$,
respectively, and both of them are of length $s$.

If $s$ is odd, then $s\ge m$ since the shortest odd
$(x_1,y_1)$-walk in $G_1$ is of length $m$;
If $s$ is even, then $s\ge n$ since the shortest even
$(x_2,y_2)$-walk in $G_2$ is of length $n$.
Therefore $d_G(x,y)=s\ge \min\{m,n\}$.
\end{pf}


\section{Main results}

Let $G$ be a connected graph with odd cycles and
${\mathscr C^o(G)}$ be the set of all odd cycles in $G$.
For $C\in {\mathscr C^o(G)}$ and $x\in V(G)$, let
 $$
 d_G(x,C)=\min\{d_G(x,y): y\in V(C)\},
 $$ and let
 $$
 \begin{array}{rl}
 & d_G^o(C)=\max\{d_G(x,C): x\in V(G-C)\} \ \ {\rm for}\ C\in {\mathscr C^o(G)},\\
 & l^o(G)=\min\{2d_G^o(C)+|V(C)|-1: C\in {\mathscr C^o(G)}\}.
 \end{array}
 $$
We define $l^o(G)=\infty$ if $G$ is bipartite.

\begin{thm}\label{thm3.1}
$\gamma(G)\le l^o(G)$ for any connected graph $G$.
\end{thm}

\begin{pf}
If $G$ contains no odd cycles, then $l^o(G)=\infty$, and so the conclusion holds.
Suppose that $G$ contains odd cycles. By Lemma~\ref{lem2.1}, $G$ is primitive.
We only need to prove that for any two vertices $x$ and $y$
in $G$, $\gamma(G;x,y)\le l^o(G)$.

By definition, there exists an odd cycle $C$ such that
 $
 l^o(G)=2d_G^o(C)+|V(C)|-1.
 $
Let $d_1=d_G(x,C)$ and $d_2=d_G(y,C)$.
Then $d_1\le d_G^o(C)$ and $d_2\le d_G^o(C)$.
Let $P_x=(x,x_1,\ldots,x_{d_1})$
and $P_y=(y,y_1,\ldots,y_{d_2})$ be two shortest paths from
$x$ and $y$ to $C$, respectively, where $x_{d_1},y_{d_2}\in V(C)$ (maybe $x_{d_1}=y_{d_2}$).
Two vertices $x_{d_1}$ and $y_{d_2}$ partition $C$ into two paths $P_1$
and $P_2$ with lengths $p_1$ and $p_2$, respectively.
Then $p_1$ and $p_2$ have different parity, say $p_1\ge p_2$.
Thus, $P_x\cup P_1 \cup P_y$ and $P_x\cup P_2 \cup P_y$ are
two $(x,y)$-walks with length of different parity and at most
 $$
  d_1+d_2+p_1\le 2d_G^o(C)+ |V(C)|=l^o(G)+1.
 $$
By Lemma~\ref{lem2.2}, $\gamma(G;x,y)\le l^o(G)$.
\end{pf}

\begin{cor}\label{cor2.10}
If $G$ is a connected graph with loops and diameter $d$,
then $\gamma(G)\le 2d$.
\end{cor}

Let $H_{n,p}$ and $F_{n,p}$ ($p\ge 1$) be two graphs,
which are obtained by joining a complete graph
$K_p$ and a cycle $C_p$ to the end-vertex $x_{n-p}$ of a path
$P_{n-p}=(x_1,x_2,\ldots,x_{n-p})$ with an edge, respectively.

The following result can be deduced by
Theorem~\ref{thm3.1}.

\begin{cor}\label{cor3.1}\textnormal{(Wang and Wang~\cite{ww93})}
Let $G$ be a primitive graph with order $n$ and odd girth $p\ge 3$.
Then $\gamma(G)\le 2n-p-1$ with equality if and only if
$G$ is isomorphic to $F_{n,p}$.
\end{cor}

\begin{pf}
By Lemma~\ref{lem2.1}, $G$ is connected and contains an odd cycle $C$ with
$l^o(G)=2d_G^o(C)+|V(C)|-1$. Since $d_G^o(C)\le n-|V(C)|$
and $|V(C)|\ge p$, by Theorem~\ref{thm3.1},
we have that
 \begin{equation}\label{e3.1}
 \begin{array}{rl}
 \gamma(G)&\le l^o(G)=2d_G^o(C)+|V(C)|-1\\
 &\le 2(n-|V(C)|)+|V(C)|-1\\
 &\le 2n-p-1.
 \end{array}
 \end{equation}
The equality implies that all equalities
in (\ref{e3.1}) hold, in particular, $d_G^o(C)=n-|V(C)|$
and $|V(C)|=p$. Thus, there is a vertex $x_1$ such that
$d_G(x_1,C)=n-p$ in $G$. Suppose $P=(x_1,x_2,\ldots,x_{n-p},x_{n-p+1})$
is a shortest path from $x_1$ to $C$, where $x_{n-p+1}$ is in $C$.
By the minimality of $P$ and primitivity of $G$, it is easy to see
that $G$ is isomorphic
to $F_{n,p}$. Also, if $G$ is isomorphic to $F_{n,p}$, then
the shortest odd closed $(x_1,x_1)$-walk is of length
$2(n-p)+p=2n-p$. This implies there is no closed $(x_1,x_1)$-walk of length
$2n-p-2$. Hence, $\gamma(F_{n,p})\ge 2n-p-1$.
\end{pf}

\begin{cor}\label{cor3.2}
If $p\ge 3$, then $\gamma(H_{n,p})=2n-2p+2$.
\end{cor}

\begin{pf}
Let $G=H_{n,p}$. Since $G$ contains $K_p$ and $p\ge 3$,
$G$ is primitive by Lemma~\ref{lem2.1}, and so
$d_G^o(C)=d_G(x_1,C)=n-p$ for any $C\in {\mathscr C^o(G)}$.
Let $C$ be a cycle of length $3$ in $G$.
By Theorem~\ref{thm3.1},
 $$
 \gamma(G)\le l^o(G)\le 2d_G^o(C)+|V(C)|-1=2(n-p)+|V(C)|-1\le 2n-2p+2.
 $$

It is clear that the shortest odd closed $(x_1,x_1)$-walk
is of length $2(n-p)+3$. This implies there is no closed
$(x_1,x_1)$-walk of length $2(n-p)+1$. Hence,
$\gamma(G)\ge 2n-2p+2$.
The conclusion follows.
\end{pf}

\begin{thm}\label{thm3.2}
Let $G_i$ be a connected graph with diameter
$d_i\ge 1$ and exponent $\gamma_i=\gamma(G_i)$ for $i=1,2$,
$G_1$ contains odd cycles, and $G=G_1 \otimes G_2$.
Then the diameter $d(G)$ of $G$ satisfies the following properties.

(1) $d(G)\ge \max\{d_1,d_2\}$.

(2) If $G_2$ contains odd cycles, then
 $$
 d(G)\ge \left\{
 \begin{array}{ll}
 \gamma_1\ & {\rm if}\ \gamma_1=\gamma_2;\\
 \min\{\gamma_1,\gamma_2\}+1\ & {\rm if}\ \gamma_1\ne \gamma_2.
 \end{array}\right.
 $$

(3) $d(G)\le \max\{\gamma_1,\gamma_2\}$.

(4) $d(G)\le \min\{\max\{\gamma_1+1,d_2\}, \max\{\gamma_2+1,d_1\}\}$
with equality if $G_2$ is bipartite.
\end{thm}

\begin{pf}
Since both $G_1$ and $G_2$ are connected and
$G_1$ contains odd cycles, by Lemma~\ref{lem2.1} and Lemma~\ref{lem2.5},
$\gamma_1$ is well-defined and $G$ is connected.
Since $d_1\ge 1$ and $d_2\ge 1$, the order of $G_1$ and
$G_2$ are no less than 2.

(1) For $i=1,2$, let $x_i$ and $y_i$ be two vertices
in $G_i$ with $d_{G_i}(x_i,y_i)=d_i$ and let
$P=(x_1x_2,\ldots,u_1u_2,\ldots,y_1y_2)$ be a shortest $(x_1x_2,y_1y_2)$-path
in $G$. Then $(x_1,\ldots,u_1,\ldots,y_1)$ and
$(x_2,\ldots,u_2,\ldots,y_2)$ are two walks in $G_1$ and $G_2$,
respectively. Thus $d(G)\ge d(P)\ge \max\{d_1,d_2\}$.

(2) Since $G_2$ contains odd cycles, $\gamma_2$ is well-defined by
Lemma~\ref{lem2.1}. Without loss of generality, assume
$\gamma_2\ge \gamma_1$ and $\gamma_1$ is odd.
By Lemma~\ref{lem2.6}, there exist
two different vertices $x_1$ and $y_1$ such that the shortest
even $(x_1,y_1)$-walk is of length $\gamma_1+1$ in $G_1$; also
there exist
two vertices $x_2$ and $y_2$ such that the shortest
odd $(x_2,y_2)$-walk is of length $\gamma_2$ or $\gamma_2+1$ in $G_2$.
By Lemma~\ref{lem2.7},
$d_G(x_1x_2,y_1y_2)\ge \min\{\gamma_1+1,\gamma_2\}$,
and so
$$
 d(G)\ge \left\{
 \begin{array}{ll}
 \gamma_1\ & {\rm if}\ \gamma_1=\gamma_2;\\
 \min\{\gamma_1,\gamma_2\}+1\ & {\rm if}\ \gamma_1\ne \gamma_2.
 \end{array}\right.
 $$

(3) 
Without loss of generality, suppose that $\gamma_2$ is well-defined
and $\gamma_2\le \gamma_1$.
Let $x=x_1x_2$ and $y=y_1y_2$ be any two different vertices in $G$.
By definition of $\gamma$,
there exist an $(x_1,y_1)$-walk and an $(x_2,y_2)$-walk of length
$\gamma_1$ in $G_1$ and $G_2$, respectively.
By Lemma~\ref{lem2.5}, there exists an $(x,y)$-walk of length $\gamma_1$,
and hence $d(G;x,y)\le \gamma_1$. By
the arbitrariness of $x$ and $y$, we have $d(G)\le \gamma_1$.

(4) 
Without loss of generality, suppose that $\gamma_2$ is well-defined, and
only need to prove $d(G)\le \max\{\gamma_1+1,d_2\}$.
Let $x=x_1x_2$ and $y=y_1y_2$ be any two different vertices
in $G$ and $d_2'=d_{G_2}(x_2,y_2)$ (maybe $x_2=y_2$).
If $d_2'\ge \gamma_1$, then there exists an $(x_1,y_1)$-walk of length
$d_2'$ in $G_1$ by definition of $\gamma$. By Lemma~\ref{lem2.5},
there exists an $(x,y)$-walk of length $d_2'$ in $G$.
If $d_2'<\gamma_1$, then one of $d_2'+\gamma_1$ and $d_2'+\gamma_1+1$
is even. By definition of $\gamma$, there exist two $(x_1,y_1)$-walks
of lengths $\gamma_1$ and $\gamma_1+1$ in $G_1$, respectively.
By Lemma~\ref{lem2.5}, there exists an $(x,y)$-walk of length
no more than $\gamma_1+1$ in $G$. Thus $d_G(x,y)\le \max\{\gamma_1+1,d_2\}$,
and hence $d(G)\le \max\{\gamma_1+1,d_2\}$ by arbitrariness of $x$ and $y$.

Now assume that $G_2$ is bipartite.
Let $x_2$ and $y_2$ be two vertices in different parts in $G_2$.
Then any $(x_2,y_2)$-walk and closed $(x_2,x_2)$-walk are of odd
and even length in $G_2$, respectively.
If $\gamma_1=1$, then $d(G)\ge d_G(x_1x_2,y_1x_2)\ge 2=\gamma_1+1$ for any
two different vertices $x_1$ and $y_1$ in $G_1$ since $|V(G_1)|\ge 2$.
Next, assume $\gamma_1\ge 2$.

By using the Lemma~\ref{lem2.6}, we have the following conclusions.
If $\gamma_1$ is odd, then there exist two different vertices $x_1$ and $y_1$
such that the shortest even $(x_1,y_1)$-walk is of length $\gamma_1+1$ in $G_1$,
and hence $d(G)\ge d_G(x_1x_2,y_1x_2)\ge \gamma_1+1$.
If $\gamma_1$ is even, then there exist two vertices $x_1$ and $y_1$
such that the shortest odd $(x_1,y_1)$-walk is of length $\gamma_1+1$ in $G$,
and hence $d(G)\ge d_G(x_1x_2,y_1y_2)\ge \gamma_1+1$.
By the conclusion (1), $d(G)\ge d_2$, and hence $d(G)=\max\{\gamma_1+1,d_2\}$.

The theorem follows.\end{pf}

\begin{cor}
Let $G_i$ be a connected graph with diameter $d_i\ge 1$ and
$l_i=l^o(G_i)$ for $i=1,2$, $G=G_1 \otimes G_2$. Then
$$d(G)\le \min\{\max\{l_1+1,d_2\}, \max\{l_2+1,d_1\}\}.$$
\end{cor}

\begin{pf}
Without loss of generality, we can suppose that both $G_1$ and $G_2$ contain odd cycles.
By Theorem~\ref{thm3.1}, $\gamma(G_1)\le l_1$
and $\gamma(G_2)\le l_2$. The conclusion follows by
the conclusion (4) in Theorem~\ref{thm3.2}.
\end{pf}

\begin{cor}
Let $G_i$ be a connected graph with diameter $d_i\ge 1$ for $i=1,2$
and $G=G_1 \otimes G_2$.
If $G_1$ contains odd cycles,
then $d(G)\le \max\{2d_1+1,d_2\}$.
\end{cor}

\begin{pf}
By Lemma~\ref{lem2.3}, $\gamma(G_1)\le 2d_1$.
The Theorem follows by the conclusion (4) in Theorem~\ref{thm3.2}.
\end{pf}\\

The following result, obtained by Leskovec et al.~\cite{lckfc10},
can be deduced by Theorem~\ref{thm3.2} immediately.

\begin{cor}\textnormal{(Leskovec et al.~\cite{lckfc10})}
Let $G_i$ be a connected graph with diameter $d_i\ge 1$
and there is a loop on every vertex of $G_i$ for $i=1,2$.
Then $d(G_1\otimes G_2)=\max\{d_1,d_2\}$.
\end{cor}

\begin{pf}
It is clear that $\gamma(G_1)=d_1$ and $\gamma(G_2)=d_2$
since each of $G_1$ and $G_2$ has a loop on every vertex.
The conclusion follows by the conclusions (1) and (3) in Theorem~\ref{thm3.2}.
\end{pf}

\begin{cor}
Let $G$ be a primitive graph with order $n\ge 2$.
Then $\gamma(G)=d(G\otimes K_2)-1$.
\end{cor}

\vspace{10pt}
By Theorem~\ref{thm3.2}, we immediately obtain our main results in this paper.

\begin{thm}\label{thm3.3}
Let $G_i$ be a connected graph with diameter $d_i\ge 1$ and exponent
$\gamma_i=\gamma(G_i)$ for $i=1,2$. If $G_1$ contains odd cycles,
then
 $$
 d(G_1\otimes G_2)=\left\{ \begin{array}{ll}
\gamma_1\ & {\rm if }\ \gamma_1=\gamma_2; \\
\max\{\gamma_2+1,d_1\}\ & {\rm if }\ \gamma_1>\gamma_2;\\
\max\{\gamma_1+1,d_2\}\ & {\rm if }\ \gamma_1<\gamma_2.
\end{array}
 \right.
 $$
\end{thm}

In Theorem~\ref{thm3.3}, we consider the diameter of the Kronecker
product of two graphs $G_1$ and $G_2$ with order no less than 2.
Next, we consider the case that at least one of $G_1$ and $G_2$ with
order 1. Let $G$ be a connected graph with order $n$ and no parallel
edges. We have noted in Section 2, $\gamma(G)=1$ if and only if
$G\cong K_n^+$. For a graph $H$ with order 1, if $G\otimes H$ is
connected, then $H\cong K_1^+$ since $G\otimes K_1$ is empty. It is
easy to see that $K_1^+\otimes G\cong G$ and then $d(K_1^+\otimes
G)=d(G)$.

\vspace{10pt}
In the following, we show the diameters for some special
Kronecker product of two graphs only by using the diameters
of factor graphs.

\begin{thm}\label{thm3.4}
Let $G_i$ be a connected graph with order $n_i\ge 2$ for $i=1,2$. Then
$d(G_1\otimes G_2)=1$ if and only if $G_1\cong K_{n_1}^+$ and $G_2\cong K_{n_2}^+$.
\end{thm}

\begin{pf}
The sufficiency is obviously.

Now we show the necessity. By contradiction.
Without loss of generality, assume $G_1\ncong K_{n_1}^+$.
Then either there exists a vertex $x$ such that it does not contain
a loop or $d(G_1)\ge 2$. Then $d(G)\ge d_G(xy,xz)\ge 2$
for any two different vertices $y,z\in V(G_2)$ or $d(G)\ge d(G_1)\ge 2$
by the conclusion (1) in Theorem~\ref{thm3.2}.
\end{pf}

\begin{thm}\label{thm3.5}
Let $G\ncong K_n^+$ be a connected graph with order $n\ge 2$ and $m\ge 2$.
Then
$$d(K_m^+\otimes G)=\left\{ \begin{array}{ll}
2, \ \ \ \ \  &  d(G)=1;\\
d(G),   &  d(G)\ge 2.
\end{array}
 \right.$$
\end{thm}
\begin{pf}
The Theorem follows by Theorem~\ref{thm3.3} since
$\gamma(K_m^+)=1$ and $\gamma(G)\ge 2$.
\end{pf}

\begin{thm}
Let $G$ be a connected graph with diameter $d\ge 1$ and
$H$ be a complete $t$ partite graph with $t\ge 3$.
Then
$$d(G\otimes H)=\left\{ \begin{array}{ll}
d,\  &  d\ge 3; \\
2,  &  d\le 2\ {\rm and}\ \gamma(G)\le 2;\\
3, &  d\le 2\ {\rm and}\ \gamma(G)>2.
\end{array}
 \right. $$
\end{thm}
\begin{pf}
It is clear that $d(H)\ge 2$, $H$ is primitive and $\gamma(H)=2$.
The Theorem follows by Theorem~\ref{thm3.3}.
\end{pf}

\begin{cor}
Let $G$ be $H_{n,p}$ or $F_{n,p}$ with odd cycles and diameter
$d_1\ge 1$, and $H$ be any connected graph with diameter $d_2\ge 1$.

(1) If $H$ is bipartite, then $d(G\otimes H)=\max\{2d_1+1,d_2\}$.

(2) If $H=G_{n_2,p_2}$ is non-bipartite, then
 $$
 d(G\otimes H)=\left\{ \begin{array}{ll}
 2d_1\ & {\rm if}\  d_1=d_2; \\
\max\{d_1,2d_2+1\}\ & {\rm if}\ d_1>d_2;\\
\max\{d_2,2d_1+1\}\ & {\rm if}\ d_1<d_2.
\end{array}
 \right. $$
\end{cor}
\begin{pf}
By Lemma~\ref{lem2.1}, $G$ is primitive since $G$ contains odd cycles.
By Lemma~\ref{lem2.5}, $G\otimes H$ is connected.
By  Corollary~\ref{cor3.1} and~\ref{cor3.2}, $\gamma(G)=2d_1$.
If $H$ is bipartite, then $H$ is not primitive by Lemma~\ref{lem2.1}.
Thus $\infty=\gamma(H)>\gamma(G)$, and hence,
$d(G\otimes H)=\max\{2d_1+1,d_2\}$ by Theorem~\ref{thm3.3}.
If $H=G_{n_2,p_2}$ is non-bipartite, then $\gamma(H)=2d_2$.
The conclusion follows by Theorem~\ref{thm3.3} immediately.
\end{pf}

\begin{cor}
Let $C_m$ be an odd cycle and $H$ be a connected graph with
order $n$ and diameter $d\ge 1$.

(1) If $H$ is bipartite, then $d(C_m\otimes H)=\max\{m,d\}$.
Hence $d(C_m\otimes P_n)=\max\{m,n-1\}$, and
$d(C_m\otimes C_n)=\max\{m,\frac{n}{2}\}$ if $n$ is even.

(2) If $H=C_n$ and $n$ is odd, then
 $$
 d(C_m\otimes C_n)=\left\{ \begin{array}{ll}
m-1\ & {\rm if}\ m=n, \\
\max\{n,\frac{m-1}{2}\}\ & {\rm if}\ m>n\\
\max\{m,\frac{n-1}{2}\}\ & {\rm if}\ m<n.
\end{array}
 \right. $$
\end{cor}


\end{document}